\documentclass[11pt]{amsart}
\usepackage{amssymb, amsthm, amscd,
enumerate, amsmath}

\newtheorem{thm}{Theorem}[section]
\newtheorem{cor}[thm]{Corollary}

\theoremstyle{definition}
\newtheorem{defn}[thm]{Definition}

\theoremstyle{remark}
\newtheorem{rem}[thm]{Remark}

\numberwithin{equation}{section}

\newcommand{\ep}{\varepsilon}

\newcommand{\la}{\lambda}

\newcommand{\si}{\sigma}
\newcommand{\Si}{\Sigma}

\newcommand{\csi}{\xi}

\newcommand{\x}{\times}

\newcommand{\CC}{\mathcal C}

\renewcommand{\SS}{{\mathfrak {F}}}

\newcommand{\imm}{{\mathrm {Imm}}}

\newcommand{\R}{\mathbb R}

\newcommand{\del}{\partial}

\newcommand{\co}{\colon\thinspace}

\begin{document}
\mathsurround=1pt 
\title[Fold maps, stably framed manifolds and immersions]
{Cobordism of fold maps, stably framed manifolds and immersions}



\thanks{The author is supported by the JSPS Research Fellowship for Young Scientists.}

\author{Boldizs\'ar Kalm\'{a}r}

\address{Kyushu University, Faculty of Mathematics, 6-10-1 Hakozaki, Higashi-ku, Fukuoka 812-8581, Japan}
\email{kalmbold@yahoo.com}


\begin{abstract}
We give complete geometric invariants of cobordisms of fold maps with oriented singular set
and cobordisms of even codimensional fold maps.
These invariants are given in terms of cobordisms of stably framed manifolds and
cobordisms of immersions with prescribed normal bundles. 
\end{abstract}

\maketitle

\section*{Introduction}

We defined and used geometric cobordism invariants of fold maps in \cite{Kal4, Kal6, Kal7}
which describe the immersion of the singular set of the fold map into the target manifold 
together with more detailed informations about
the tubular neighbourhood of the singular set of the fold map in the source manifold.
In this paper we define further invariants 
which describe the cobordism class of the source manifold and its fold map outside 
of the singular set as well. 
In \cite{Kal4} we showed that the cobordism groups of fold maps contain
stable homotopy groups of spheres as direct components and in \cite{Kal7} we showed
that they also contain stable homotopy groups of the classifying spaces $BO(k)$. In this paper, by
using the results of Ando \cite[Theorem~0.1, Theorem~3.2]{An3} about the existence of 
fold maps, we show that together with our invariants defined in \cite{Kal4, Kal6, Kal7},
 cobordism groups of manifolds with stable framings (see Section~\ref{ujcob}) give complete cobordism
invariants of 
the cobordism classes of even codimensional fold maps and
the cobordism classes of framed fold maps (see Definition~\ref{framedcobdef}) 
of arbitrary codimension into the Euclidean space (see Theorem~\ref{invarithm},
Corollary~\ref{foldcobcor}).

We emphasize that our results in this paper are based on
elementary constructions of geometric cobordism invariants \cite{Kal4, Kal6, Kal7},
constructions of special fold maps \cite{Kal7} and an elementary application (see Section~\ref{completebiz})
of a simple corollary \cite[Theorem~0.1, Theorem~3.2]{An3} of the 
h-principle for fold maps of Ando \cite[Theorem~0.5, Theorem~2.1]{An3}.

Independently from our present paper Sadykov \cite{Sad4} gave a splitting of the cobordism groups of fold maps in terms of homotopy groups of spectra, also in the odd codimension case. These spectra were constructed \cite{Sad2} using the
h-principle for fold maps of Ando \cite[Theorem~0.5, Theorem~2.1]{An3} but with a quite different and much more sophisticated approach than ours. Sadykov also showed that our geometric invariants defined in \cite{Kal4, Kal7, Kal6} coincide with direct summands of the splitting given by \cite{Sad4}.

For other more sophisticated 
applications of the h-principle for fold maps of Ando 
toward cobordisms of fold maps, see, \cite{An2, An5, An4, An, An6}.

For further results about cobordisms of fold maps with positive codimension, see,
\cite{EkSzuTer, Szucs4}. For cobordisms of negative codimensional 
singular maps with completely different
approach from our present paper, see, \cite{An2, An5, An4, An, An6,
Ik, IS, Kal, Kal2, Sad3, Saspecgen, Sa}.

The paper is organized as follows.
In Section~\ref{preli} we give the basic definitions,
in Section~\ref{mainres} we state our main results and
in Section~\ref{completebiz} we prove our results about complete invariants.

The author would like to thank Prof. Yoshifumi Ando for the interesting and useful conversations.

\subsection*{Notations}
In this paper the symbol ``$\amalg$'' denotes the disjoint union, 
for any number $x$ the symbol ``$\lfloor x \rfloor$'' denotes the greatest
integer $i$ such that $i \leq x$,
$\ep^1_X$ denotes the trivial line bundle over the space $X$,
$\ep^1$ denotes the trivial line bundle over the point,
and the symbols $\csi^k$, $\eta^k$, etc. usually denote $k$-dimensional real vector bundles.
The symbol $T\csi^k$ denotes
the Thom space of the bundle $\csi^k$.
The symbol $\imm^{\csi^k}_{N}(n-k,k)$ denotes
the cobordism group of $k$-codimensional immersions into 
an $n$-dimensional manifold $N$
whose normal bundles are induced from $\csi^k$ (this
group is isomorphic to the group $\{\dot N, T\csi^k \}$, where
$\dot N$ denotes the one point compactification of the manifold $N$
and the symbol $\{X,Y \}$ denotes the group of stable homotopy classes of continuous
maps from the space $X$ to the space $Y$). 
The symbol $\imm^{\csi^k}(n-k,k)$ denotes
the cobordism group of $k$-codimensional immersions into $\R^n$ 
whose normal bundles are induced from $\csi^k$ (this
group is isomorphic to $\pi_{n}^s(T\csi^k)$).
The symbol $\pi_n^s(X)$ ($\pi_n^s$) denotes the $n$th stable homotopy group of the space $X$ (resp. spheres).
The symbols $\Omega_{m}$ and ${\mathfrak {N}}_{m}$ denote 
the usual cobordism groups of oriented and unoriented $m$-dimensional manifolds, respectively.
The symbol ``id$_A$'' denotes the identity map of the space $A$.
The symbol $\ep$ denotes a small positive number.
All manifolds and maps are smooth of class $C^{\infty}$.

\section{Preliminaries}\label{preli}

\subsection{Fold maps}
 
Let $n \geq 1$ and $q \geq 0$.
Let $Q^{n+q}$ and $N^n$ be smooth manifolds of dimensions $n+q$ and $n$ 
respectively. Let $p \in Q^{n+q}$ be a singular point of 
a smooth map $f \co Q^{n+q} \to N^{n}$. The smooth map $f$  has a {\it fold 
singularity of index $\la$} at the singular point $p$ if we can write $f$ in some local coordinates around $p$  
and $f(p)$ in the form 
\[  
f(x_1,\ldots,x_{n+q})=(x_1,\ldots,x_{n-1}, -x_n^2 - \cdots -x_{n+\la-1}^2 + x_{n+\la}^2 + \cdots + x_{n+q}^2)
\] 
for some $\la$ $(0 \leq \la \leq q+1)$ (the index $\la$ is well-defined if
we consider that $\la$ and $q+1-\la$ represent the same index). 

A smooth map $f \co Q^{n+q} \to N^{n}$ is called a {\it fold map} if $f$ has only 
fold singularities.

A smooth map $f \co Q^{n+q} \to N^n$ 
  has a {\it definite fold
singularity} at a fold singularity $p \in Q^{n+q}$ if $\la = 0$ or $\la = q+1$,
otherwise $f$ has an {\it indefinite fold singularity of index $\la$}
at the fold singularity $p \in Q^{n+q}$.

Let $S_{\la}(f)$ denote the set of fold singularities of index $\la$ of $f$ in $Q^{n+q}$.
Note that $S_{\la}(f) = S_{q+1-\la}(f)$.
 Let $S_f$ denote the set $\bigcup_{\la} S_{\la}(f)$.

Note that the set $S_f$ is an ${(n-1)}$-dimensional submanifold of the manifold
$Q^{n+q}$.

Note that each connected component of the manifold $S_f$ has its own index $\la$ if
we consider that $\la$ and $q+1-\la$ represent the same index. 

Note that for a fold map $f \co Q^{n+q} \to N^{n}$ and for an index $\la$ ($0 \leq \la \leq \lfloor q/2 \rfloor$)
the codimension one immersion $f |_{S_{\la}(f)} \co S_{\la}(f) \to N^n$ 
of the singular set of index $\la$ $S_{\la}(f)$ has a canonical framing 
(i.e., trivialization of the normal bundle) by identifying canonically the set of 
fold singularities of index $\la$ 
of the map $f$ with   
the fold germ 
$(x_1,\ldots,x_{n+q}) \mapsto (x_1,\ldots,x_{n-1}, -x_n^2 - \cdots -x_{n+\la-1}^2 + x_{n+\la}^2 + \cdots + x_{n+q}^2)$,
see, for example, \cite{Sa1}.

If $f \co Q^{n+q} \to N^n$ is a fold map in general position, then 
the map $f$
restricted to the singular set $S_f$ is a general positional
 codimension one immersion  into the target manifold $N^n$.
 
Since every fold map is in general position after a small perturbation, 
and we study maps under the equivalence relation {\it cobordism}
(see Definition~\ref{cobdef}),
in this paper we can restrict ourselves to studying fold maps which are 
in general position.
Without mentioning we suppose that a fold map $f$ is in general position.

\begin{defn}[Framed fold map]\label{framedfoldmap}
We say that a fold map $f \co Q^{n+q} \to N^n$ is {\it framed} if
the codimension one immersion $f|_{S_f} \co S_f \to N^n$ 
of the singular set $S_f$ has a framing 
(i.e., trivialization of the normal bundle) such that
for each index $\la$ with $0 \leq \la \leq \lfloor q/2 \rfloor$
the framing of $f|_{S_{\la}(f)} \co S_{\la}(f) \to N^n$
coincides with the canonical framing.
For odd $q$,
the framing of the immersion of the singular set of index $(q+1)/2$
can be arbitrary.
\end{defn}

Note that if we have a framed fold map $f$ into an oriented manifold $N^n$, then 
the singular set $S_f$ has a natural orientation which gives at any point of $f(S_f)$ (together
with the framing of the immersion of $S_f$) the orientation of the target
manifold $N^n$. 

\begin{defn}[Oriented fold map]\label{orifoldmap}
A fold map $f \co Q^{n+q} \to N^{n}$ is {\it oriented} if there is a chosen consistent orientation
of every fiber at their regular points. 
\end{defn}

For example, a  fold map $f \co Q^{n+q} \to N^{n}$ between oriented manifolds is canonically oriented.

Note that an oriented fold map $f \co Q^{n+q} \to N^{n}$ with odd $q$ may not have a framing in the sense of Definition~\ref{framedfoldmap}.

%

\subsection{Equivalence relations of fold maps}\label{kob}

\begin{defn}[Cobordism]\label{cobdef} 
Two fold maps $f_i \co Q_i^{n+q} \to N^n$ $(i=0,1)$  
of closed $({n+q})$-dimensional manifolds $Q_i^{n+q}$ $(i=0,1)$ 
into an $n$-dimensional manifold $N^n$ are  
{\it cobordant} if 
\begin{enumerate}[a)]
\item
there exists a fold map 
$F \co X^{n+q+1} \to N^n \times [0,1]$ of a compact $(n+q+1)$-dimensional 
manifold $X^{n+q+1}$, 
\item
$\del X^{n+q+1} = Q_0^{n+q} \amalg Q_1^{n+q}$ and 
\item
${F|}_{Q_0^{n+q} \x [0,\ep)}=f_0 \x
{\mathrm {id}}_{[0,\ep)}$ and ${F|}_{Q_1^{n+q} \x (1-\ep,1]}=f_1 \x 
{\mathrm {id}}_{(1-\ep,1]}$, where 
$Q_0^{n+q} \x [0,\ep)$
 and $Q_1^{n+q} \x (1-\ep,1]$ are small collar neighbourhoods of $\del X^{n+q+1}$ with the
identifications $Q_0^{n+q} = Q_0^{n+q} \x \{0\}$ and $Q_1^{n+q} = Q_1^{n+q} \x \{1\}$. 
\end{enumerate}
If the fold maps $f_i \co Q_i^{n+q} \to N^n$ $(i=0,1)$  are oriented,
we say that they are {\it oriented cobordant} (or shortly {\it cobordant} if it is clear from the context) if
they are cobordant in the above sense via an oriented fold map $F$, such that the orientations are compatible on the boundary.

We call the map $F$ a {\it cobordism} between $f_0$ and $f_1$.
\end{defn}

This clearly defines an equivalence relation on the set of (oriented) fold maps 
of closed $({n+q})$-dimensional manifolds into an  
$n$-dimensional manifold $N^n$.

We denote 
 the set of (oriented) cobordism classes of (oriented) fold maps of closed $({n+q})$-dimensional manifolds 
into an $n$-dimensional manifold 
$N^n$
by $\CC ob_{N,f}^{}(n+q,-q)$ (resp. $\CC ob_{N,f}^{O}(n+q,-q)$). When $N^n=\R^n$, we denote it by
$\CC ob_{f}^{}(n+q,-q)$ (resp. $\CC ob_{f}^{O}(n+q,-q)$).
We note that we can define a commutative semigroup operation in the usual way on the 
set of cobordism classes $\CC ob_{N,f}^{}(n+q,-q)$ (resp. $\CC ob_{N,f}^{O}(n+q,-q)$)
by the disjoint union.
If the target manifold $N^n$ is the Euclidean space $\R^n$ (or more
generally if $N^n$ has the form $\R^1 \x M^{n-1}$ for some $(n-1)$-dimensional manifold  
$M^{n-1}$), then the elements in the semigroup $\CC ob_{N,f}^{}(n+q,-q)$
(resp. $\CC ob_{N,f}^{O}(n+q,-q)$)
have their inverses: namely, compose them with a reflection
in a hyperplane (in $\{0\} \x M^{n-1}$ in general, see \cite{Szucs4}).
Hence the semigroups $\CC ob_{N,f}(n+q,-q)$ (resp. $\CC ob_{N,f}^{O}(n+q,-q)$) are in this case 
actually groups.

\begin{defn}[Framed cobordism]\label{framedcobdef}
Two (oriented) framed fold maps $$f_i \co Q_i^{n+q} \to N^n$$ $(i=0,1)$  
of closed $({n+q})$-dimensional manifolds $Q_i^{n+q}$ $(i=0,1)$ 
into an $n$-dimensional manifold $N^n$ are  
{\it (oriented) framed cobordant} if 
they are (oriented) cobordant in the sense of Definition~\ref{cobdef}
by a framed fold map $F \co X^{n+q+1} \to N^n \times [0,1]$ 
such that the framing of the immersion 
$F|_{S_F} \co S_F \to N^n \times [0,1]$ restricted 
to the immersion $f_i|_{S_{f_i}} \co S_{f_i} \to N^n \times \{ i \}$
coincides with the framing of the fold map $f_i$ $(i=0,1)$.
\end{defn}

Let us denote the (oriented) framed cobordism semigroup of (oriented)
framed fold maps of $(n+q)$-dimensional manifolds
into $N^n$ by $\CC ob_{N,f,fr}^{}(n+q,-q)$ (resp. $\CC ob_{N,f,fr}^{O}(n+q,-q)$).

Note that for even codimension $q = 2k$ ($k \geq 0$)
the (oriented) framed cobordism semigroup 
$\CC ob_{N,f,fr}^{}(n+q,-q)$ (resp. $\CC ob_{N,f,fr}^{O}(n+q,-q)$) is naturally isomorphic to
the (oriented) cobordism semigroup
$\CC ob_{N,f}^{}(n+q,-q)$ (resp. $\CC ob_{N,f}^{O}(n+q,-q)$).

\subsection{Cobordism invariants of fold maps}\label{geomcobinv}

As an imitation of lifting positive codimensional singular maps 
\cite[Section~7, Proof of Theorem~2]{Szucs4},
we defined and used geometric invariants of cobordisms of fold maps 
(for the definitions and notations, see, \cite[Section~2]{Kal7}),
namely the homomorphisms
\[
\csi_{{\la}, q}^N  \co \CC ob_{N,f}^{(O)}(n+q,-q) \to \imm^{\ep^1_{B(O(\la) \x O(q+1-\la))}}_N(n-1,1)
\]
for $0 \leq \la < (q+1)/2$ and
\[
\csi_{(q+1)/2, q}^N  \co \CC ob_{N,f}^{(O)}(n+q,-q) \to \imm^{l^1}_N(n-1,1)
\]
for $q$ odd and $\la=(q+1)/2$. These homomorphisms map
 a cobordism class of a fold map $f$ into the cobordism class of the immersion of its 
fold singular set $S_{\la}(f)$ of index $\la$  with normal bundle induced from
the target of the universal fold germ bundle of index 
$\la$. In the case of oriented fold maps (e.g., in the case 
of oriented manifolds $Q^{n+q}$ and $N^n$), we have the analogous homomorphisms 
\[
\csi_{{\la}, q}^{O,N}  \co \CC ob_{N,f}^{O}(n+q,-q) \to \imm^{\ep^1_{BS(O(\la) \x O(q+1-\la))}}_N(n-1,1)
\]
for $0 \leq \la < (q+1)/2$ and
\[
\csi_{(q+1)/2, q}^{O,N}  \co \CC ob_{N,f}^{O}(n+q,-q) \to \imm^{\tilde l^1}_N(n-1,1)
\]
for $q$ odd and $\la=(q+1)/2$ as well.
We used these homomorphisms in \cite{Kal4, Kal6, Kal7} in order to describe cobordisms of fold maps.

For $0 \leq \la \leq (q+1)/2$,
we define the (oriented) framed cobordism invariants
\[
\csi_{{\la}, q}^N  \co \CC ob_{N,f,fr}^{(O)}(n+q,-q) \to \imm^{\ep^1_{B(O(\la) \x O(q+1-\la))}}_N(n-1,1)
\]
and
\[
\csi_{{\la}, q}^{O,N}  \co \CC ob_{N,f,fr}^{O}(n+q,-q) \to \imm^{\ep^1_{BS(O(\la) \x O(q+1-\la))}}_N(n-1,1)
\]
in the analogous way.

\subsection{Framed cobordism of manifolds}\label{ujcob}

%
%
%

\begin{defn}[Stably $(n-1)$-framed manifolds and stably $(n-1)$- framed cobordism]\label{stabframecob}
For $n>0, q \geq 0$
an $(n+q)$-dimensional manifold $Q^{n+q}$
is {\it stably $(n-1)$-framed} if the vector bundle $TQ^{n+q} \oplus \ep^2_{Q^{n+q}}$
has $n+1$ sections that are linearly independent at every point of $Q^{n+q}$ (shortly, we say that it has $n+1$
independent sections).

Let $Q_i^{n+q}$ be closed (oriented) stably $(n-1)$-framed $(n+q)$-dimensional manifolds, i.e.,
the vector bundles $TQ_i^{n+q} \oplus \ep^2_{Q_i^{n+q}}$
have $n+1$ independent sections $e_i^1, \ldots, e_i^{n+1}$ ($i=0,1$). 
We say that the manifolds $Q_0^{n+q}$ and $Q_1^{n+q}$ 
are {\it stably (oriented) $(n-1)$-framed cobordant} if 
\begin{enumerate}[a)]
\item
they are (oriented) cobordant in the
usual sense by an (oriented) $(n+q+1)$-dimensional manifold $W^{n+q+1}$,
\item
one of the two trivial line bundles in the direct sum $TQ_0^{n+q} \oplus \ep^2_{Q_0^{n+q}}$ (resp. $TQ_1^{n+q} \oplus \ep^2_{Q_1^{n+q}}$)
corresponds to the inward (resp.\ outward) normal section of the boundary of
 $W^{n+q+1}$,
\item
the vector bundle $TW \oplus \ep^1_W$
has $n+1$ independent sections $f^1, \ldots, f^{n+1}$,
\item
the sections $f^j$ 
restricted to the boundary $Q_0^{n+q} \amalg Q_1^{n+q}$ of $W^{n+q+1}$ 
coincide with the sections
$e_i^j$  ($j=1, \ldots, n+1$ and $i=0,1$).
\end{enumerate}
\end{defn}

We denote the set of (oriented) stably $(n-1)$-framed cobordism classes
of closed (oriented) stably $(n-1)$-framed $(n+q)$-dimensional manifolds 
 by 
$\CC_{n+q}^{}(n)$ (resp. $\CC_{n+q}^{O}(n)$) which is an abelian group with the disjoint union as operation.

By \cite[Lemma 3.1]{Sa1} (see also \cite[Lemma 3.1]{An3}) there is a homomorphism 
$$\si^{(O)}_{n,q} \co \CC ob_{f,fr}^{(O)}(n+q,-q) \to   \CC_{n+q}^{(O)}(n),$$
which maps a cobordism class of a framed fold map $g \co Q^{n+q} \to \R^n$ 
into the stably $(n-1)$-framed cobordism class of the
source manifold $Q^{n+q}$ with the stable framing obtained by 
\cite[Lemma 3.1]{Sa1}\footnote{In \cite[Lemma 3.1]{Sa1} $n$ independent sections for 
$TQ^{n+q} \oplus \ep^1_{Q^{n+q}}$ are constructed.}.


\section{Results}\label{mainres}

In this section, the target manifold $N^n$ is always assumed to be $\R^n$.
Note that 
the group
$$\imm^{\ep^1_{B(O(\la) \x O(q+1-\la))}}_N(n-1,1)$$ (or $\imm^{\ep^1_{BS(O(\la) \x O(q+1-\la))}}_N(n-1,1)$),
which is the target of the homomorphism $\csi_{{\la}, q}^{N}$ (resp. $\csi_{{\la}, q}^{O,N}$),
is naturally identified with the group
$$\pi^s_{n-1} \oplus \pi^s_{n-1}(B(O(\la) \x O(q+1-\la)))$$ (resp. 
$\pi^s_{n-1} \oplus \pi^s_{n-1}(BS(O(\la) \x O(q+1-\la)))$).

The main results of this section are the following.

\begin{thm}\label{invarithm}
Let $n>0, q \geq 0$.
Then, the homomorphisms
\begin{multline*}
  \si^{O}_{n,q} \oplus \bigoplus_{1 \leq \la \leq (q+1)/2 } \csi_{{\la}, q}^O
\co 
\CC ob_{f, fr}^{O}(n+q,-q) \longrightarrow  \\  {\CC}_{n+q}^O{(n)}
\oplus  
\bigoplus_{1 \leq \la \leq (q+1)/2} \pi^s_{n-1} \oplus \pi^s_{n-1}(BS(O(\la) \x O(q+1-\la)))
\end{multline*}
and 
\begin{multline*}
  \si^{}_{n,q} \oplus \bigoplus_{1 \leq \la \leq (q+1)/2 } \csi_{{\la}, q}
\co 
\CC ob_{f, fr}^{}(n+q,-q) \longrightarrow  \\  {\CC}_{n+q}{(n)}
\oplus  
\bigoplus_{1 \leq \la \leq (q+1)/2} \pi^s_{n-1} \oplus \pi^s_{n-1}(B(O(\la) \x O(q+1-\la)))
\end{multline*}
denoted by $\Im_{n,q}^O$ and $\Im_{n,q}$, respectively,
are injective. In other words, for closed (oriented) manifolds
$Q_i^{n+q}$ $(i=0,1)$
two framed fold maps $$f_i \co Q_i^{n+q} \to \R^n$$ $(i=0,1)$ are (oriented) framed  
cobordant if and only if 
$$\Im_{n,q}^{(O)}([f_0]) = \Im_{n,q}^{(O)}([f_1]).$$
\end{thm}


\begin{cor}\label{foldcobcor}
For $k \geq 0$,
the homomorphism $\Im_{n,2k}^{(O)}$ gives a complete invariant
of the (oriented) cobordism group $\CC ob_{f}^{(O)}(n+2k,-2k)$ of fold maps.
\end{cor}

\begin{rem}
By Corollary~\ref{foldcobcor} the homomorphisms
$$\Im_{n,0}^{O} \co 
\CC ob_{f}^{O}(n,0) \to  {\CC}_{n}^O{(n)}$$ 
and
$$\Im_{n,0}^{} \co 
\CC ob_{f}^{}(n,0) \to  {\CC}_{n}{(n)}$$ 
are injective, and by
Ando \cite[Theorem~3.2]{An3} they are surjective as well. Hence,
we have $\CC ob_{f}^{O}(n,0) =  {\CC}_{n}^O{(n)}$ and
$\CC ob_{f}^{}(n,0) =  {\CC}_{n}{(n)}$.
Note that since the group ${\CC}_{n}^O{(n)}$ is isomorphic to $\pi^s_n$, 
we obtain another argument for the isomorphism
$\CC ob_{f}^{O}(n,0) =  \pi^s_n$ (for the original proof, see Ando \cite{An2, An}).
\end{rem}

\section{Proof}\label{completebiz}

In this section, we prove Theorem~\ref{invarithm}.

\begin{proof}[Proof of Theorem~\ref{invarithm}]
Let $f \co Q^{n+q} \to \R^n$ be a framed fold map.
If \[\csi_{{\la}, q}([f]) \in \imm^{\ep^1_{B(O(\la) \x O(q+1-\la))}}(n-1,1)\] 
is zero for $1 \leq \la \leq (q+1)/2$, then
by gluing the null-cobordisms of the fold singularity bundles \cite{Kal7} to the fold map $f$,
we obtain a framed fold map $F \co W^{n+q+1} \to \R^n \x [0,1]$,
where 
\begin{enumerate}
\item
$W^{n+q+1}$ is a compact $(n+q+1)$-dimensional manifold with boundary
$Q^{n+q} \amalg P^{n+q}$, 
\item
$F|_{Q^{n+q}\x [0,\ep)} = f \x
{\mathrm {id}}_{[0,\ep)}$, 
where 
$Q^{n+q} \x [0,\ep)$
 is a small collar neighbourhood of $Q^{n+q}$ ($\subset \del W^{n+q+1}$) in $W^{n+q+1}$ with the
identification $Q^{n+q} = Q^{n+q} \x \{0\}$ and
\item
$F$ is a submersion into $\R^n \x (0,1)$
near $P^{n+q}$.
\end{enumerate}

If $  \si^{}_{n,q}([f])$ is zero, then since
the manifolds $Q^{n+q}$ and $P^{n+q}$ are stably $(n-1)$-framed cobordant 
by the manifold $W^{n+q+1}$ and framing $\SS$ obtained by \cite[Lemma 3.1]{Sa1} from the framed fold map $F$,
the stably $(n-1)$-framed manifold $P^{n+q}$ is also zero in the 
stable $(n-1)$-framed cobordism group 
$  \CC_{n+q}^{}(n)$. Hence by gluing a stably $(n-1)$-framed null-cobordism of
$P^{n+q}$ to $W^{n+q+1}$,
we obtain a
compact $(n+q+1)$-dimensional manifold $X^{n+q+1}$ with boundary
$Q^{n+q}$ such that the bundle $TX^{n+q+1} \oplus \ep^1_{X^{n+q+1}}$ 
has an $(n+1)$-framing which
coincides with the framing $\SS$ on $W^{n+q+1}$. 

Since $\R^n \x [0,1]$ is contractible,
we can extend the map $F$ to
a continuous map $G \co X^{n+q+1} \to
\R^n \x [0,1]$. Note that the $(n+1)$-framing 
of the bundle $TX^{n+q+1} \oplus \ep^1_{X^{n+q+1}}$
gives a fiberwise epimorphism $H$ of the bundle
$TX^{n+q+1} \oplus \ep^1_{X^{n+q+1}}$ into the tangent bundle $T(\R^n \x [0,1])$ covering
the continuous map $G$
by mapping 
the $n+1$ frames to the canonical bases of $T(\R^n \x [0,1])$
at any points of $X^{n+q+1}$. 
We may suppose that near the submanifold $P^{n+q}$ of $X^{n+q+1}$
the continuous map $G$ coincides with the fold map $F$,
and by the construction of the $(n+1)$-framing $\SS$ of the bundle 
$TW^{n+q+1} \oplus \ep^1_{W^{n+q+1}}$
(see \cite[Lemma 3.1]{Sa1} and also \cite[Lemma 3.1]{An3})
the bundle homomorphism $H|_{TX^{n+q+1}} \co TX^{n+q+1} \to T(\R^n \x [0,1])$
is given by the differential of the fold map $F$ (which is a submersion)
near $P^{n+q}$.

Hence
similarly to \cite[Proof of Theorem~3.2]{An3} 
using the relative h-principle for fold maps 
\cite[Theorem~0.5, Theorem~2.1]{An3}, we see that 
there is a framed fold map $g \co X^{n+q+1} \to
\R^n \x [0,1]$ 
which coincides with $F$ on $W^{n+q+1}$ and
whose boundary\footnote{The boundary of a map $g$ from a manifold $X$ with boundary $\del X$
is the restriction $g|_{\del X}$.} coincides with the framed fold map $f$. 
Hence the framed fold map $f$ is framed
null-cobordant. 
The oriented case is proved in a similar way.
\end{proof}

\end{document}